%

\magnification=1200
\input amstex
\documentstyle{amsppt}

\NoLimitsOnNames
\define\nn{\text{$\Bbb N$} }
\define\bn{\text{$\beta\Bbb N$} }
\define\uu{{\Cal U}}
\define\vv{{\Cal V}}
\define\ww{{\Cal W}}
\define\li#1{#1\text{-}\!\lim}
\define\ct{^\frown}
\define\si{\Sigma}
\topmatter
\title
Ultrafilters: Where\\
Topological Dynamics $=$ Algebra $=$ Combinatorics
\endtitle
\rightheadtext{Ultrafilters: Dynamics $=$ Algebra $=$ 
Combinatorics}
\author
Andreas Blass
\endauthor
\address
Mathematics Dept., University of Michigan, Ann Arbor, 
MI 48109, U.S.A.
\endaddress
\email
ablass\@umich.edu
\endemail
\thanks Partially supported by NSF grant DMS-9204276.
\endthanks
\abstract
We survey some connections between topological 
dynamics, semigroups of ultrafilters, and combinatorics.  
As an application, we give a proof, based on ideas of 
Bergelson and Hindman, of the Hales-Jewett partition 
theorem.
\endabstract
\endtopmatter
\document

Furstenberg and his co-workers have shown \cite{15, 16, 
17} how to deduce combinatorial consequences from 
theorems about topological dynamics in compact metric 
spaces.  Bergelson and Hindman \cite4 applied similar 
methods in non-metrizable spaces, particularly the 
Stone-\v Cech compactification \bn of the discrete space 
of natural numbers.  This approach and related ideas of 
Carlson \cite{11} lead to particularly simple formulations 
since many of the basic concepts of dynamics, when 
applied to $\bn$, can be expressed in terms of a 
semigroup operation on $\bn$, the natural extension of 
addition on $\nn$.  The semigroup \bn can also 
substitute, in many contexts, for the enveloping 
semigroups (\cite{14}) traditionally used in topological 
dynamics.  Further simplifications and applications of 
these ideas were developed in \cite3.

The purpose of this paper is to survey some of these 
developments.  In contrast to most surveys, however, we 
include some detailed proofs, in order to emphasize their 
simplicity.  In the first three sections, we develop the 
necessary theory of dynamics and the equivalent 
semigroup structure in $\bn$.  In the fourth section, we 
apply the theory to present proofs of Hindman's partition 
theorem for finite sums and of the Hales-Jewett theorem 
about homogeneous combinatorial lines in cubes.  A final 
section (omitted for lack of time in the talk on which this 
paper is based) compares the ultrafilters discussed in the 
earlier sections with other ultrafilters traditionally 
related to combinatorics, for example selective 
ultrafilters.

\head
1. Ultrafilters
\endhead

Throughout most of this paper we are concerned with 
ultrafilters on the set \nn of natural numbers.  These are 
usually defined by some version of the following 
set-theoretic definition, in which we have included some 
redundant clauses for ease of future reference.

\definition{Definition 1}
An {\sl ultrafilter\/} on \nn is a family $\Cal U$ of 
subsets of \nn such that 
\roster
\item If $X\subseteq Y$ and $X\in\Cal U$ then $Y\in\Cal 
U$.
\item If $X,Y\in\Cal U$ then $X\cap Y\in\uu$.
\item $\emptyset\notin\uu$.
\item $\nn\in\uu$.
\item For any $X\subseteq\nn$, either $X\in\uu$ or 
$\nn-X\in\uu$.
\item If $X\cup Y\in\uu$ then either $X\in\uu$ or 
$Y\in\uu$.
\endroster
\enddefinition

The first four of these clauses define {\sl filters\/}.

The implications in \therosteritem2 and \therosteritem6 
are reversible, by \therosteritem1.

To each $a\in\nn$ is associated a {\sl principal\/} or {\sl 
trivial\/} ultrafilter, namely $\hat a= 
\{X\subseteq\nn\allowmathbreak\mid a\in X\}$.  In many contexts, we 
identify $\hat a$ with $a$. 

The following alternative definition expresses the usual 
way of viewing ultrafilters topologically.

\definition{Definition 2}
An {\sl ultrafilter\/} on \nn is a point in the Stone-\v 
Cech compactification \bn of the discrete space $\nn$.
\enddefinition

The two definitions are equivalent in the following sense.  
For any point $p\in\bn$, the family of subsets of \nn 
whose closures in \bn contain $p$ satisfies Definition 1.  
Conversely, if $\uu$ is as in Definition 1, then the 
closures in \bn of the sets in $\uu$ have exactly one 
point $p$ in common.  And the constructions described in 
the preceding two sentences are inverse to each other.

Although the two definitions of ultrafilters above are the 
most familiar ones, two other, equivalent definitions will 
be more useful for our purposes.  The first of these uses 
the notion of a {\sl quantifier\/} $Q$ over $\nn$.  This is 
an operation which applies to a formula $\varphi(n)$ 
with a free variable ranging over \nn and produces a 
new formula $(Qn)\,\varphi(n)$ in which $n$ is no longer 
free; it is required that replacing $\varphi(n)$ by an 
equivalent formula $\psi(n)$ yields an equivalent result 
$(Qn)\,\psi(n)$.  Formally, a quantifier can be identified 
with the set of those $X\subseteq\nn$ for which 
$(Qn)\,n\in X$ is true.  Under this identification, the 
following definition amounts to Definition 1.

\definition{Definition 3}
An {\sl ultrafilter\/} on \nn is a quantifier $\uu$ over 
\nn that respects the propositional connectives in the 
sense that the following equivalences hold for all 
formulas $\varphi(n)$ and $\psi(n)$
\roster
\item $(\uu n)\,\varphi(n)\land(\uu n)\,\psi(n)\iff
(\uu n)\,(\varphi(n)\land\psi(n))$
\item $(\uu n)\,\varphi(n)\lor(\uu n)\,\psi(n)\iff
(\uu n)\,(\varphi(n)\lor\psi(n))$
\item $\neg(\uu n)\,\varphi(n)\iff
(\uu n)\,\neg\varphi(n)$
\endroster
\enddefinition

If $\uu$ is an ultrafilter in the sense of Definition 1, then 
the corresponding quantifier $(\uu n)$, usually read ``for 
$\uu$-almost all $n$,'' is defined by 
$$
(\uu n)\,\varphi(n)\iff\{n\in\nn\mid\varphi(n)\}\in\uu;
$$
conversely, from a quantifier as in Definition 3 we can 
define 
$$
\uu=\{X\subseteq\nn\mid(\uu n)\,n\in X\}.  
$$

As all propositional connectives (of any number of 
arguments) can be expressed in terms of $\neg$ and 
$\land$, they are all respected by ultrafilter quantifiers.

Finally, we give another topological definition, which will 
provide a connection to dynamics.

\definition{Definition 4}
An {\sl ultrafilter\/} on \nn is a uniform operation on 
sequences in compact Hausdorff spaces.  That is, it is an 
operator assigning to every sequence $(x_n)_{n\in\nn}$ 
in every compact Hausdorff space $X$ a point 
$\li\uu_nx_n\in X$ subject to the requirement that, if 
$f:X\to Y$ is a continuous map to another compact 
Hausdorff space, then $f(\li\uu_nx_n)=\li\uu_nf(x_n)$.
\enddefinition

The easiest way to connect this definition with the 
previous ones is to notice that a sequence $(x_n)$ in $X$ 
is a (continuous) function $x$ from the discrete space \nn 
into $X$; if $X$ is a compact Hausdorff space, then this 
map extends uniquely to $\bar x:\bn\to X$, and so each 
$p\in\bn$ yields a point $\bar x(p)\in X$.  Uniformity is 
easy to check, so an ultrafilter in the sense of Definition 2 
yields one in the sense of Definition 4.  Conversely, an 
operation as in Definition 4 can be applied to the 
sequence in \bn whose $n^{\text{th}}$ term is $n$, 
yielding a point $p\in\bn$, and these constructions are 
inverse to each other.  

One can verify that $\li\uu_nx_n$ is the unique point in 
$X$ such that every neighborhood $G$ of it satisfies $(\uu 
n)\,x_n\in G$.  In particular, it follows that $\li\uu_nx_n$ 
is a limit point or a member of the sequence $(x_n)$.  
Thus, a non-trivial ultrafilter can be regarded as a 
systematic way of passing to a limit point of any 
sequence.

The trivial ultrafilter $\hat a$ corresponds in Definition 2 
to the point  $a\in\nn\subseteq\bn$, in Definition 3 to 
the ``quantifier'' that just substitutes $a$ for the 
quantified variable, and in Definition 4 to the operation 
that picks out the $a^{\text{th}}$ term from sequences.

The set \bn of all ultrafilters on \nn admits a binary 
operation $+$, extending ordinary addition on \nn (see 
for example \cite{12, 23}).  In the context of Definition 4, 
it amounts to an iteration of limit operations:
$$
\li{(\uu+\vv)}_px_p=\li\uu_m(\li\vv_nx_{m+n}).
$$
Translating this into the language of quantifiers, one 
again finds an iteration:
$$
((\uu+\vv)p)\,\varphi(p)\iff
(\uu m)(\vv n)\,\varphi(m+n).
$$
The equivalent characterizations in terms of Definitions 1 
and 2 are more complicated, at least on first sight.  
Definition 1 leads to
$$
\uu+\vv=\{X\subseteq\nn\mid\{m\mid\{n\mid
m+n\in X\}\in\vv\}\in\uu\}.
$$
And for Definition 2 we have the following description of 
addition.  Start with ordinary addition 
$+:\nn\times\nn\to\nn$.  Extend it by continuity to 
$+:\nn\times\bn\to\bn$, fixing the first argument in \nn 
and requiring continuity in the second.  Then, fixing the 
second argument in \bn and requiring continuity in the 
first, obtain an extension $+:\bn\times\bn\to\bn$.

Notice that the operation $+$ on \bn is a continuous 
function of the left summand for any fixed value in \bn 
of the right summand, but it is not a continuous function 
of the right summand for a fixed left summand unless the 
latter is in $\nn$ (see \cite{23, Section 10}).  I refer to 
continuity in the left argument as left-continuity, and I 
therefore call \bn a left topological semigroup.  (Caution: 
Some authors use ``right''  instead of ``left'' because the 
right translations are continuous, and some authors 
define $\uu+\vv$ to be what I would call $\vv+\uu$; 
authors who disagree with me on both points therefore 
say ``left,'' just as I do, though they mean the opposite.)

The addition operation on \bn is associative (most easily 
checked using the quantifier description of $+$); it is 
commutative as long as one of the summands is in $\nn$, 
but not in general (for details, see \cite{23, Section 10}).

\head
2. Dynamics
\endhead

Topological dynamics is concerned with the behavior of 
iterations of a continuous map $T$ from a space $X$ into 
itself.  (Actually, it is considerably more general \cite{14}, 
but the preceding description covers what will be 
relevant here.)  For the purposes of this paper, a {\sl 
dynamical system\/} consists of a compact Hausdorff 
space $X$ and a continuous function $T:X\to X$.  We write 
$T^n$ for the $n^{\text{th}}$ iterate $T\circ 
T\circ\dots\circ T$ of $T$.  To study the limiting behavior 
of these iterates for large $n$, we define (as in \cite{24}) 
for each ultrafilter $\uu$ on \nn
$$
T^\uu(x)=\li\uu_nT^n(x).
$$
Regarded as a function of $\uu\in\bn$, for fixed $x\in X$, 
this is the continuous extension to \bn of the function 
$\nn\to X:n\mapsto T^n(x)$.  It follows that 
$\{T^\uu(x)\mid \uu\in\bn\}$ is the closure in $X$ of the 
forward orbit $\{T^n(x)\mid n\in\nn\}$ of the point $x$.  

But regarded as a function of $x$ for fixed $\uu$, 
$T^\uu(x)$  need not be continuous unless $\uu$ is 
principal.  If $\uu$ is the principal ultrafilter $\hat n$, 
then $T^\uu=T^n$, so no confusion will be caused by 
identifying $\hat n$ with $n$ in this context.  

The notion of iteration with respect to an ultrafilter, 
$T^\uu$, connects nicely with the addition operation on 
ultrafilters in that 
$$
T^\uu(T^\vv(x))=T^{\uu+\vv}(x).
$$
Indeed, we have
$$\align
T^{\uu+\vv}(x)&=\li{(\uu+\vv)}_pT^p(x)\\
&=\li\uu_m\li\vv_nT^{m+n}(x)= 
\li\uu_m\li\vv_nT^m(T^n(x))\\
&=\li\uu_mT^m(\li\vv_nT^n(x))
\qquad\text{(as $T^m$ is continuous)}\\
&=T^\uu(T^\vv(x)).
\endalign$$

We next introduce some concepts from topological 
dynamics, i.e., concepts about the behavior of $T^n(x)$ for 
large $n$.  In each definition, it is assumed that $(X,T)$ is 
a dynamical system.  More information about these 
concepts can be found in \cite{14, 15}.

\definition{Definition}
A point $x\in X$ is {\sl recurrent\/} if, for each  
neighborhood $G$ of $x$, infinitely many $n\in\nn$ 
satisfy $T^n(x)\in G$.  It is {\sl uniformly recurrent\/} if, 
for each neighborhood $G$ of $x$, there is $M\in\nn$ 
such that $\forall n\,\exists k<M\,T^{n+k}(x)\in G$.
\enddefinition

Thus, recurrence means that, under the iteration of $T$, 
the point $x$ returns to each of its neighborhoods 
infinitely often.  Uniform recurrence bounds how long the 
sequence of iterates can stay out of any given 
neighborhood; there is $M$ depending on $G$ such that of 
every $M$ consecutive iterates at least one is in $G$.

\definition{Definition}
Two points $x,y\in X$ are {\sl proximal\/} if, for every 
neighborhood $G$ of the diagonal in $X\times X$, 
infinitely many $n\in\nn$ satisfy $(T^n(x),T^n(y))\in G$.
\enddefinition

Proximality is usually defined in the context of metric 
spaces by requiring that, for every positive 
$\varepsilon$, infinitely many $n$ have the distance 
between $T^n(x)$ and $T^n(y)$ smaller than 
$\varepsilon$.  This definition clearly makes use not of 
the full metric structure but only of the associated 
uniform structure; that is, it makes sense in any uniform 
space.  A compact Hausdorff space has a unique uniform 
structure, and the definition we gave for compact 
Hausdorff spaces is just the specialization to this case of 
the general concept in uniform spaces.

The dynamical concepts just defined can be elegantly 
expressed in terms of ultrafilter iterations, as follows.

\proclaim{Theorem  1}
Let $(X,T)$ be a dynamical system.
\roster
\item 
A point $x\in X$ is recurrent if and only if $T^\uu(x)=x$ 
for some non-trivial ultrafilter $\uu$ on $\nn$, if and only 
if $T^\uu(x)=x$ for some $\uu\neq\hat0$.
\item 
A point $x\in X$ is uniformly recurrent if and only if for 
every ultrafilter $\vv$ on \nn there is an ultrafilter 
$\uu$ on \nn with $T^\uu(T^\vv(x))=x$
\item
Two points $x,y\in X$ are proximal if and only if there is 
an ultrafilter $\uu$ on \nn with $T^\uu(x)=T^\uu(y)$, if 
and only if there is a non-trivial ultrafilter $\uu$ on \nn 
with $T^\uu(x)=T^\uu(y)$.
\endroster
\endproclaim

\demo{Proof}
\therosteritem1\quad
By definition, recurrence means that $x$ is a limit point 
of the sequence $(T^n(x))$.  That implies that $x$ is in the 
closure $\{T^\uu(x)\mid\uu\in\bn-\{\hat0\}\}$ of 
$\{T^n(x)\mid n\in\nn-\{0\}\}$.  This in turn implies that 
$x\in\{T^\uu(x)\mid\uu\in\bn-\nn\}$, i.e., that we can 
take $\uu$ non-trivial.  Indeed, if $\uu$ were trivial, say 
$\uu=\hat n\neq\hat0$, then $T^n(x)=x$, so $T^{nk}(x)=x$ 
for all $k$, and therefore, if we take $\uu'$ to be a 
non-principal ultrafilter containing the set of multiples of 
$n$, then $T^{\uu'}(x)=x$ also.  Finally, 
$x\in\{T^\uu(x)\mid\uu\in\bn-\nn\}$ implies 
recurrence, since every $T^\uu(x)$ with non-principal 
$\uu$ is a limit point of $\{T^n(x)\mid n\in\nn\}$.

\therosteritem2\quad
Assume first that $x$ is uniformly recurrent, and let an 
ultrafilter $\vv$ be given.  Temporarily fix a closed 
neighborhood $G$ of $x$, and let $M$ be as in the 
definition of uniform recurrence for this neighborhood.  
So $\forall n\,\exists k<M\,T^{n+k}(x)\in G$.  As only 
finitely many $k$'s occur and as $\vv$ is an ultrafilter, 
the same $k$ must work for $\vv$-almost all $n$.  Fix 
this $k$, so we have $(\vv n)\,T^{n+k}(x)\in G$.  
Equivalently, $(\vv n)\,T^n(x)\in T^{-k}(G)$.  As $T^{-
k}(G)$ is closed, $T^\vv(x)\in T^{-k}(G)$, and so $T^k 
(T^\vv(x))\in G$. Now un-fix $G$, and remember that, as 
$X$ is a compact Hausdorff space, every neighborhood of 
$x$ includes a closed neighborhood.  So we have shown 
that, for every neighborhood $G$ of $x$, the set 
$$
Y_G=\{k\in\nn\mid T^k(T^\vv(x))\in G\}
$$
is non-empty.  Clearly, $Y_{G_1}\cap Y_{G_2}=Y_{G_1\cap 
G_2}$, so, as $G$ ranges over the neighborhoods of $x$, 
the sets $Y_G$ generate a filter.  Extend it to an ultrafilter 
$\uu$.  Then we have, for each neighborhood $G$ of $x$, 
$(\uu k)\,T^k(T^\vv(x))\in G$, so $T^\uu(T^\vv(x))=x$, as 
desired.

Conversely, suppose $x$ is not uniformly recurrent, and 
fix an open neighborhood $G$ such that no $M$ satisfies 
the definition of uniform recurrence.  That is, for all 
$M\in\nn$, the set
$$
Y_M=\{n\in\nn\mid(\forall k<M)\,T^{n+k}(x)\notin G\}
$$
is nonempty.  These sets $Y_M$ generate a filter, as they 
form a chain, so there is an ultrafilter $\vv$ containing 
all of them.  For every $k\in\nn$ we have, since 
$Y_{k+1}\in\vv$, 
$(\vv n)\,T^{n+k}(x)\notin G$; so $(\vv n)\,T^n(x)\notin 
T^{-k}(G)$; so, as $T^{-k}(G)$ is open, $T^\vv(x)\notin T^{-
k}(G)$; so $T^k(T^\vv(x))\notin G$.  As $k$ is arbitrary, it 
follows that, for any ultrafilter $\uu$, 
$T^\uu(T^\vv(x))\notin G$, as desired.

\therosteritem3\quad
Suppose $T^\uu(x)=T^\uu(y)$.  If $\uu$ is principal, then 
proximality follows trivially, so suppose $\uu$ is 
non-principal.  To prove that $x$ and $y$ are proximal, 
let $G$ be any neighborhood of the diagonal.  Since 
$(T^\uu(x),T^\uu(y))=\li\uu_n(T^n(x),T^n(y))$ is on the 
diagonal by assumption, there must be infinitely many 
$n$ (in fact $\uu$-almost all $n$) such that 
$(T^n(x),T^n(y))\in G$, as required.

Conversely, suppose $x$ and $y$ are proximal.  So, as $G$ 
ranges over neighborhoods of the diagonal, the sets 
$$
Y_G=\{n\in\nn\mid(T^n(x),T^n(y))\in G\}
$$
are non-empty, and they generate a filter (because 
$Y_{G_1}\cap Y_{G_2}=Y_{G_1\cap G_2}$), which we 
extend to an ultrafilter $\uu$.  Then we have, for all 
closed neighborhoods $G$ of the diagonal, $(\uu n)\, 
(T^n(x),T^n(y))\in G$ and so $(T^\uu(x),T^\uu(y))\in G$.  
But the intersection of all closed neighborhoods of the 
diagonal is just the diagonal, so we conclude that 
$T^\uu(x)=T^\uu(y)$.
\qed\enddemo

We close this section by pointing out a simpler connection 
between dynamical systems and ultrafilters: ultrafilters 
provide the universal example of a dynamical system.  
The compact Hausdorff space \bn with the shift map 
$S:\bn\to\bn:\uu\mapsto\hat1+\uu$ is a dynamical 
system and enjoys the following universal property.  If 
$(X,T)$ is any dynamical system and $x$ is any point in 
$X$, then there is a unique continuous map $f:\bn\to X$ 
such that $f\circ S=T\circ f$ and $f(0)=x$, namely the 
map defined by $f(\uu)=T^\uu(x)$.  Thus, $(\bn,S)$ may 
be regarded as the free dynamical system on one 
generator.  

Henceforth, when we refer to \bn as a dynamical system, 
we mean $(\bn,S)$.

\head
3. Dynamics $=$ Algebra
\endhead

The addition operation defined for \bn in Section ~1 is 
just the ultrafilter iteration of the shift map, i.e., of the 
universal dynamical system.  Indeed, we have, for any 
ultrafilters $\uu$ and $\vv$ on $\nn$
$$
S^\uu(\vv)=\li\uu_nS^n(\vv)=\li\uu_n(n+\vv)=\uu+\vv,
$$
where at the last step we used that $\li\uu$ commutes 
with continuous maps, like addition as a function of its 
left argument, and that $\li\uu_nn=\uu$.

Of course, this equivalence between iteration in $(\bn,S)$ 
and addition allows us to reformulate Theorem~1 
algebraically.  We do so in the following theorem, adding 
some more reformulations in terms of subsemigroups 
(i.e., non-empty subsets closed under addition) and ideals 
in the semigroup $(\bn,+)$.  A {\sl left ideal\/} is a 
non-empty set $I\subseteq\bn$ such that if $\uu\in I$ 
and $\vv\in\bn$ then $\vv+\uu\in I$; {\sl right ideals\/} 
and {\sl two-sided ideals\/} are defined analogously.

\proclaim{Theorem 2}
\roster
\item
An ultrafilter $\uu$ is recurrent in \bn if and only if 
$\vv+\uu=\uu$ for some $\vv\neq\hat0$.
\item
An ultrafilter $\uu$ is uniformly recurrent in \bn if and 
only if, for each $\vv$, there is $\ww$ with 
$\ww+\vv+\uu=\uu$, if and only if $\uu$ belongs to a 
minimal (closed) left ideal in $\bn$.
\item
Two ultrafilters $\uu_1$ and $\uu_2$ are proximal in 
\bn if and only if there is an ultrafilter $\vv$ such that 
$\vv+\uu_1=\vv+\uu_2$.
\item
An ultrafilter $\uu$ generates a minimal closed 
subsemigroup of \bn if and only if it is idempotent, i.e., 
$\uu+\uu=\uu$.
\endroster
\endproclaim

\demo{Proof}
Parts \therosteritem1, \therosteritem3, and the first part 
of \therosteritem2 are immediate consequences of the 
corresponding parts of Theorem~1 and the fact that 
$S^\uu(\vv)=\uu+\vv$.  

To finish the proof of \therosteritem2, notice first that 
every ultrafilter $\uu$ generates a left ideal, namely 
$\bn+\uu=\{\vv+\uu\mid\vv\in\bn\}$.  It follows that a 
minimal left ideal, being the ideal generated by any of its 
elements $\uu$, is closed, for it is the image of the 
compact space \bn under the continuous map adding 
$\uu$ on the right.  That is why ``closed'' is parenthesized 
in \therosteritem2; putting it in or leaving it out doesn't 
affect the statement.  Now to say that a particular 
ultrafilter $\uu$ is in a minimal left ideal is to say that 
the left ideal $\bn+\uu$ that it generates is minimal or, 
equivalently, is generated by each of its elements.  That 
is, for every ultrafilter $\vv$, the ideal $\bn+\vv+\uu$ 
generated by $\vv+\uu$ must be all of $\bn+\uu$.  
Equivalently,  $\bn+\vv+\uu$ must contain the generator 
$\uu$ of $\bn+\uu$.  But that means that, for every 
$\vv$, we can express $\uu$ as $\ww+\vv+\uu$ by 
suitably choosing $\ww$.  This completes the proof of 
\therosteritem2.

\therosteritem4, which is included in the theorem 
because of its analogy to \therosteritem2, is trivial in one 
direction, as $\{\uu\}$ is a closed subsemigroup if $\uu$ 
is idempotent.  To prove the non-trivial direction (due, as 
far as I know, to Ellis \cite{14}), let $C$ be a minimal 
closed subsemigroup of \bn and let $\uu\in C$.  Then 
$C+\uu$ is also closed (being the image of the compact set 
$C$ under a continuous map) and a subsemigroup of $C$, 
so by minimality it equals $C$.  In particular it contains 
$\uu$.  So the set $D=\{\vv\in C \mid\vv+\uu=\uu\}$ is 
nonempty.  It is closed (being the pre-image of $\{\uu\}$ 
under a continuous map) and also a subsemigroup of $C$, 
so it equals $C$ and therefore contains $\uu$.  That is, 
$\uu+\uu=\uu$, and the proof is complete.  (It follows, of 
course, by minimality, that $C=\{\uu\}$.)
\qed\enddemo

The information in Theorem~2 about ideals and 
subsemigroups can be used to give quick existence proofs 
for the corresponding sorts of ultrafilters.

\proclaim{Corollary  }
There exist uniformly recurrent ultrafilters.  There exist 
non-trivial idempotent ultrafilters.
\endproclaim

\demo{Proof}
The intersection of a chain of closed subsemigroups of 
\bn is again a closed subsemigroup; it is non-empty by 
compactness, and it is obviously closed topologically and 
closed under addition.  By Zorn's Lemma, there are 
minimal closed subsemigroups of $\bn-\nn$.  Their 
elements are idempotent by \therosteritem4 of the 
theorem.  The same argument applied to closed left ideals 
yields uniformly recurrent points.
\qed\enddemo

The concepts characterized in Theorem~2 are related to 
each other as follows.

\proclaim{Theorem  3}
Each of the following statements about an ultrafilter 
$\uu$ and the dynamical system \bn implies the next.
\roster
\item
$\uu$ is uniformly recurrent and proximal to 0.
\item
$\uu$ is idempotent.
\item
$\uu$ is recurrent and proximal to 0.
\endroster
\endproclaim

\demo{Proof}
\therosteritem1$\implies$\therosteritem2\quad
Let $\uu$ be uniformly recurrent and proximal to 0.  By 
Theorem~2\therosteritem3, fix $\vv$ with 
$\vv+\uu=\vv+0=\vv$.  By Theorem~2\therosteritem2, 
fix $\ww$ with 
$$
\ww+\vv+\uu=\uu.
$$
Combining these two equations, we get $\uu=\ww+\vv$, 
and substituting this into the displayed equation we get 
$\uu+\uu=\uu$.

\therosteritem2$\implies$\therosteritem3\quad
If $\uu$ is idempotent, then the requirement 
$\vv+\uu=\uu$ for recurrence and the requirement 
$\vv+\uu=\vv$ for proximality to 0 (see Theorem~2) are 
satisfied by taking $\vv=\uu$.
\qed\enddemo

The preceding results connect the algebraic properties of 
\bn with its dynamical properties, but in fact, thanks to 
the universality of \bn among dynamical systems, we 
easily get connections between the algebra of \bn and 
arbitrary dynamical systems.

\proclaim{Theorem  4}
Let $(X,T)$ be a dynamical system and let $x\in X$.  If 
$\uu$ is (uniformly) recurrent in \bn then $T^\uu(x)$ is 
(uniformly) recurrent in $X$.  If\/ $\uu_1$ and $\uu_2$ 
are proximal in \bn then $T^{\uu_1}(x)$ and 
$T^{\uu_2}(x)$ are proximal in $X$.
\endproclaim

\demo{Proof}
Each part is proved by combining the corresponding parts 
of Theorems~1 and 2 with the fact that $T^{\uu+\vv}= 
T^\uu\circ T^\vv$.

Suppose $\uu$ is recurrent in $\bn$.  So by 
Theorem~2\therosteritem1 there is $\vv$ with 
$\vv+\uu=\uu$.  Then 
$T^\vv(T^\uu(x))=T^{\vv+\uu}(x)=T^\uu(x)$, so $T^\uu(x)$ 
is recurrent by Theorem~1\therosteritem1.

Suppose $\uu$ is uniformly recurrent.  By 
Theorem~2\therosteritem2, for every $\vv$ there is 
$\ww$ with $\ww+\vv+\uu=\uu$ and therefore 
$T^\ww(T^\vv(T^\uu(x)))=T^\uu(x)$.  By 
Theorem~1\therosteritem2, $T^\uu(x)$ is uniformly 
recurrent.

Finally, suppose $\uu_1$ and $\uu_2$ are proximal.  By 
Theorem~2\therosteritem3, there is $\vv$ with 
$\vv+\uu_1=\vv+\uu_2$.  Then 
$T^\vv(T^{\uu_1}(x))=T^\vv(T^{\uu_2}(x))$.  By 
Theorem~1\therosteritem3, $T^{\uu_1}(x)$ and 
$T^{\uu_2}(x)$ are proximal.
\qed\enddemo

As an application of these connections between dynamics 
and algebra, we give a short proof of the Auslander-Ellis 
Theorem \cite{15}.

\proclaim{Theorem  5}
Let $(X,T)$ be a dynamical system.  For each $x\in X$, 
there exists a uniformly recurrent $y$ proximal to $x$.
\endproclaim

\demo{Proof}
By the corollary of Theorem~2, there exists a uniformly 
recurrent $\vv\in\bn$.  It follows immediately that 
every ultrafilter of the form $\ww+\vv$ is uniformly 
recurrent.  The set $\bn+\vv$ of such ultrafilters is a 
closed subsemigroup of $\bn$.  By Zorn's Lemma, it 
includes a minimal closed subsemigroup.  By 
Theorem~2\therosteritem2, there is an idempotent 
$\uu\in\bn+\vv$.  Then $\uu$, being uniformly 
recurrent and idempotent, is also proximal to 0 by 
Theorem~3.  

Now for $X$, $T$, and $x$ as in the theorem, let 
$y=T^\uu(x)$.  Then, by Theorem~4, $y$ is uniformly 
recurrent and proximal to $T^0(x)=x$.
\qed\enddemo

The property of ultrafilters, ``uniformly recurrent and 
proximal to 0,'' which played a key role in the proof of 
Theorem 5, has alternative algebraic descriptions that 
will be useful later.  To introduce them, we first define a 
partial ordering of the idempotent ultrafilters by 
$$
\uu\leq\vv\iff\uu+\vv=\vv+\uu=\uu.
$$
This definition and parts of the next theorem are from 
\cite4.  When we refer to an idempotent ultrafilter as {\sl 
minimal\/}, we mean with respect to this ordering.

\proclaim{Theorem  6}
The following three assertions are equivalent, for any 
$\uu\in\bn$.
\roster
\item
$\uu$ is uniformly recurrent and proximal to 0.
\item
$\uu$ is idempotent and belongs to some minimal left 
ideal of $\bn$.
\item
$\uu$ is a minimal idempotent.
\endroster
Furthermore, these equivalent conditions imply that 
$\uu$ belongs to every two-sided ideal of $\bn$.  Finally, 
every idempotent ultrafilter $\uu$ is $\geq$ a minimal 
idempotent.
\endproclaim

\demo{Proof}

The equivalence of \therosteritem1 and \therosteritem2 
is immediate from Theorem~2\therosteritem2 and 
Theorem~3.

To prove \therosteritem2$\implies$\therosteritem3, 
assume \therosteritem2, and suppose $\vv$ is an 
idempotent $\leq\uu$.  Since $\uu$ is uniformly 
recurrent by Theorem~2, choose $\ww$ so that 
$\ww+\vv+\uu=\uu$, which reduces, in view of 
$\vv\leq\uu$, to $\ww+\vv=\uu$.  Using this, the 
idempotence of $\vv$, and again $\vv\leq\uu$, we 
compute
$$
\vv=\uu+\vv=\ww+\vv+\vv=\ww+\vv=\uu,
$$
so $\uu$ is minimal.

We next show that, if $\uu$ is idempotent and 
$I\subseteq\bn+\uu$ is a minimal left ideal, then there is 
an idempotent $\vv\leq\uu$ in $I$.  Since we already 
know that \therosteritem2 implies \therosteritem3, this 
gives the last sentence of the theorem; it will also be 
useful in establishing 
\therosteritem3$\implies$\therosteritem2. So let such 
$\uu$ and $I$ be given.  Being a closed subsemigroup of 
$\bn$, $I$ contains an idempotent $\ww$.  Being in 
$\bn+\uu$, this $\ww$ satisfies $\ww+\uu=\ww$ 
because $\uu$ is idempotent.  Let $\vv=\uu+\ww$.  Then 
$\vv$ belongs to the left ideal $I$ because $\ww$ does.  
From $\ww+\uu=\ww$ and the idempotence of $\ww$ 
and $\uu$, we infer
$$
\vv+\uu=\uu+\ww+\uu=\uu+\ww=\vv,
$$
$$
\uu+\vv=\uu+\uu+\ww=\uu+\ww=\vv,
$$
and
$$
\vv+\vv=\uu+\ww+\uu+\ww=\uu+\ww+\ww=
\uu+\ww=\vv,
$$
which mean that $\vv\leq\uu$, as desired.

The proof of \therosteritem3$\implies$\therosteritem2 is 
now easy.  If $\uu$ is a minimal idempotent, apply Zorn's 
Lemma to get a minimal left ideal $I\subseteq\bn+\uu$ 
as in the preceding paragraph, and let $\vv$ be obtained 
as there.  Being $\leq\uu$, this $\vv$ must be equal to 
$\uu$ by minimality.  So $\uu\in I$.  

Finally, we must prove that every ultrafilter satisfying 
\therosteritem2 belongs to every two-sided ideal.  In 
fact, that every minimal left ideal $I$ is included in every 
two-sided ideal $J$.  To see this, let $\uu\in I$ and 
$\vv\in J$.  Then $I\cap J$ is non-empty because it 
contains $\vv+\uu$.  So $I\cap J$ is a left ideal, and it 
must equal $I$ because $I$ is minimal.  So $I\subseteq 
J$.
\qed\enddemo

\head
4. Combinatorics
\endhead

In this section, we apply the results obtained above to 
give relatively easy proofs of some highly non-trivial 
combinatorial theorems.  The first of these is Hindman's 
Theorem, first proved in \cite{22}.  A simpler proof was 
given by Baumgartner \cite{1}, but we shall give two yet 
simpler (given the preceding machinery) arguments, one 
due to Furstenberg \cite{15} and the other to Galvin and 
Glazer \cite{12, 19, 23}.

\proclaim{Theorem  7}
If \nn is partitioned into finitely many pieces, then there 
is an infinite $H\subseteq\nn$ such that all finite sums of 
distinct members of $H$ lie in the same piece.
\endproclaim

\demo{Furstenberg's Proof}
Let the given partition have $k$ pieces, and regard it as a 
function $\nn\to K$, where $K$ is a $k$-element set.  Let 
$X$ be the set of all functions $\nn\to K$, topologized by 
giving $K$ the discrete topology and then giving $X$ the 
product topology.  Thus, $K$ is a compact Hausdorff 
space, and the given partition is a point $x\in X$.  Let 
$T:X\to X$ be the shift map, defined by $T(y)(n)=y(n+1)$; 
it is clearly continuous, so we have a dynamical system.  
By Theorem~5, let $y\in X$ be uniformly recurrent and 
proximal to $x$.  We write out what these properties of 
$y$ mean for our specific $X$ and $T$.  Uniform 
recurrence means that, given any $n\in\nn$, there is 
$N\in\nn$ such that the initial segment 
$(y(0),y(1),\dots,y(n-1))$ of $y$ recurs at least once in 
every segment $(y(r),\dots,y(r+N-1))$ of $y$ of length 
$N$.  Proximality means that, given any $N$, there are 
infinitely many intervals of length $N$ where $x$ and 
$y$ agree, $(x(r),\dots,x(r+N-1))=(y(r),\dots,y(r+N-1))$.

Let $c=y(0)$.  We intend to complete the proof by finding 
infinitely many natural numbers, all of whose finite sums 
are mapped to $c$ by $x$.

By uniform recurrence, find $N_0$ such that $c$ occurs at 
least once among every $N_0$ consecutive terms in $y$.  
By proximality, find a place, beyond term 0, where $N_0$ 
consecutive terms of $x$ coincide with those of $y$ and 
therefore contain a $c$.  So we can fix $h_0>0$ with 
$x(h_0)=c$.  This $h_0$ will be the first member of our 
$H$.

By uniform recurrence, find $N_1$ such that among 
every $N_1$ consecutive terms of $y$ there are $h_0+1$ 
consecutive terms that coincide with $y(0),\dots,y(h_0)$; 
in particular, among every $N_1$ consecutive terms, 
there are two terms a distance $h_0$ apart where $y$ 
has the value $c$ (the same as at 0 and $h_0$).  By 
proximality, there are two places a distance $h_0$ apart 
where $x$ has value $c$, say $x(h_1)=x(h_1+h_0)=c$, 
with $h_1>h_0$.  $h_1$ will be the next member of $H$.

Repeating this process, we inductively choose $h_i$ so 
that, for all sums $s$ of zero or more elements of 
$\{h_0,\dots,h_{i-1}\}$, we have $x(s+h_i)=y(s+h_i)=c$.  
This is done by finding $N$ such that every $N$ 
consecutive terms of $y$ contain a segment that coincides 
with the initial segment of $y$ up to the largest $s$, and 
then finding a segment of length $N$ beyond $h_{i-1}$  
where $x$ and $y$ coincide.

The set of all finite sums of distinct $h_n$'s is clearly 
included in the partition piece corresponding to $c$.
\qed\enddemo

\demo{Galvin's and Glazer's Proof}
This proof uses ultrafilters directly rather than applying 
them via Theorem 5 to dynamics on other spaces.  Let 
$\uu$ be any idempotent non-trivial ultrafilter on \nn 
(by the corollary to Theorem~2), and let $C$ be the piece 
of the partition that is in $\uu$ (by clause 
\therosteritem6 of Definition~1).  So we have $(\uu 
n)\,n\in C$.  As $\uu$ is idempotent, we also have (cf. the 
quantifier form of the definition of $+$ in $\bn$) $(\uu 
n)\,(\uu k)\,n+k\in C$.  As ultrafilter quantifiers respect 
propositional connectives,
$$
(\uu n)\,[n\in C \land(\uu k)\,n+k\in C].
$$
So we can fix $h_0$ with $h_0\in C$ and (re-naming 
variables) $(\uu n)\,h_0+n\in C$.  Using again that $\uu$ 
is idempotent and respects connectives, we find 
$$
(\uu n)\,[n\in C\land
(\uu k)\,n+k\in C\land
h_0+n\in C\land
(\uu k)\,h_0+n+k\in C].
$$
So we can fix $h_1$ having the four properties listed for 
$n$ inside the brackets.  In particular, $h_1\in C$ and 
$h_0+h_1\in C$.  

Repeating this process, we inductively choose $h_i$ so 
that, for all sums $s$ of zero or more elements of 
$\{h_0,\dots,h_{i-1}\}$, we have $s+h_i\in C$ and $(\uu 
n)\, s+h_i+n\in C$.  The latter property, when expanded 
by idempotence, ensures that it is possible to choose 
$h_{i+1}$ to keep the induction going. (In fact, 
$\uu$-almost all numbers can serve as $h_{i+1}$.)  
Clearly,
all finite sums of distinct members of $H=\{h_n\mid 
n\in\nn\}$ are in $C$.
\qed\enddemo

Notice that the Galvin-Glazer proof shows that the piece 
that contains the homogeneous $H$ can be taken to be 
any piece of the given partition that belongs to some 
idempotent ultrafilter.

We turn next to an application of these ideas in the 
context of words over a finite alphabet, rather than 
natural numbers.  We shall prove the Hales-Jewett 
Theorem \cite{20, 19}, but first we need some definitions 
and notational conventions.  

Let $\si$ be a finite set, which we call an alphabet, and 
let $W$ be the set of words on $\si$, i.e., the set of finite 
sequences of members of $\si$.  Let $v$ be an object, 
called a variable, that is not in $\si$; let $A$ be the set of 
words on $\si\cup\{v\}$; and let $V=A-W$.  So $V$ is the 
set of words on $\si\cup\{v\}$ in which $v$ actually 
occurs; these are often called {\sl variable words\/} over 
$\si$.  For each $a\in\si$, we define a function $\hat a: 
A\to W$ sending each $x\in A$ to the result of 
substituting $a$ for $v$ in $x$; we call $\hat a(x)$ an {\sl 
instance\/} of $x$.  Notice that if $x\in W$ then $\hat 
a(x)=x$.  

With this notation, the Hales-Jewett Theorem \cite{20} is 
as follows.

\proclaim{Theorem  8}
If $W$ is partitioned into finitely many pieces, then there 
is an $x\in V$ whose instances all lie in the same piece of 
the partition.
\endproclaim

The proof is a special case of arguments from \cite3.

\demo{Proof}
First, observe that $A$ is a semigroup under the 
operation ${}\ct$ of concatenation, that $W$ is a 
subsemigroup, that $V$ is a two-sided ideal in $A$, and 
that each $\hat a$ is a homomorphism $A\to W$.  The 
operation ${}\ct$ can be extended to the Stone-\v Cech 
compactification $\beta A$ just as addition was extended 
to $\bn$.
$$
(\uu\ct\vv x)\,\varphi(x)\iff
(\uu y)\,(\vv z)\,\varphi(y\ct z).
$$
It is easy to verify that, in the compact left-topological 
semigroup $\beta A$, $\beta W$ is a closed 
subsemigroup, $\beta V$ is a closed, two-sided ideal, and 
the continuous extension of $\hat a$, which we still call 
$\hat a:\beta A\to\beta W$, is a homomorphism.

The algebraic results about \bn proved earlier generalize 
easily to semigroups like $\beta A$ and $\beta W$.  We 
apply the analogs in this context of several parts of 
Theorem~6.  In particular, there is a minimal idempotent 
$\ww\in\beta W$.  In $\beta A$, this $\ww$ is 
idempotent but not necessarily minimal.  (In fact, we 
shall see in  a moment that it is definitely not minimal.)  
There is a minimal idempotent $\vv\leq\ww$ in $\beta 
A$.  It belongs to every two-sided ideal, so $\vv\in\beta 
V$.  (In particular, $\vv\neq\ww$.)  

For any $a\in\si$, since $\hat a$ is a homomorphism 
$\beta A\to\beta W$, we can infer from $\vv\leq\ww$ 
that $\hat a(\vv)\leq\hat a(\ww)=\ww$ (the last 
equality because $\hat a$ is the identity on $W$ and 
hence on $\beta W$).  By minimality of $\ww$ in $\beta 
W$, it follows that $\hat a(\vv)=\ww$.

Now let $W$ be partitioned into finitely many pieces, and 
let $X$ be the piece that is in the ultrafilter $\ww$.  For 
each $a\in\si$, we have $\hat a^{-1}(X)\in\vv$ because 
$X\in\ww=\hat a(\vv)$.  So $\bigcap_{a\in\si}\hat a^{-
1}(X)$ is non-empty.  (In fact it is in $\vv$.)  Any 
element $x$ of this intersection clearly serves as the $x$ 
required in the theorem.
\qed\enddemo

The proof actually establishes a stronger theorem, obtained by
broadening the notion of ``instance'' to allow a specified, finite set
of words in $W$ (not merely single letters) as the $a$'s being
substituted for $v$.  Unlike Theorem 8, this stronger form is
non-trivial even in the case where $\Sigma$ consists of just one
letter; indeed this case amounts to van der Waerden's theorem on
arithmetic progressions.  (Van der Waerden's theorem is usually
deduced from Theorem 8 for a $b$-element alphabet by using base $b$
expansions of natural numbers; see
\cite{19}.)

\head
5. P-Points
\endhead

In this section, we briefly discuss the connections 
between the ultrafilters discussed earlier and other, 
perhaps more familiar (from \cite{2, 3, 6, 8, 10, 13, 25, 26} 
for example), special ultrafilters.  We begin with a pair of 
definitions.

\definition{Definition}
A non-principal ultrafilter $\uu$ on $\nn$ is {\sl 
selective\/} if every function on $\nn$ becomes 
one-to-one or constant when restricted to a suitable set 
in $\uu$.
\enddefinition

\definition{Definition}
A non-principal ultrafilter $\uu$ on $\nn$ is a {\sl 
P-point\/} if every function on $\nn$ becomes 
finite-to-one or constant when restricted to a suitable set 
in $\uu$.
\enddefinition

Both of these definitions have numerous equivalent 
forms.  The versions above were chosen to make it 
obvious that all selective ultrafilters are P-points.  The 
definition of P-point is just a combinatorial reformulation 
of the usual topological notion of P-point specialized to 
the space $\bn$: a point such that the intersection of any 
countably many neighborhoods is again a (not necessarily 
open) neighborhood.

Selectivity also has a topological formulation, based on a theorem of
Kunen \cite{10} that characterizes selective ultrafilters $\uu$ as
those enjoying the following Ramsey property.  If the set $[\nn]^2$ of
two-element subsets of $\nn$ is partitioned into two pieces, then
there is a set $H\in\uu$ all of whose two-element subsets lie in one
piece.  (The corresponding statement for $[\nn]^k$ holds for all
finite $k$, and there are even infinitary generalizations; see
\cite{26}.)  If we identify $[\nn]^2$ with the ``above diagonal''
subset $\{(a,b)\mid a<b\}$ of $\nn^2$, then this Ramsey property says
that $[\nn]^2$ together with the sets $H\times H$ for $H\in\uu$
generate an ultrafilter on $\nn\times\nn$.  Let
$\tau:\beta(\nn\times\nn)\to\bn\times\bn$ be the continuous extension
of the inclusion map $\nn\times\nn\to\bn\times\bn$.  Then ultrafilters
on $\nn\times\nn$ that contain $H\times H$ for all $H\in\uu$ are
precisely those sent by $\tau$ to $(\uu,\uu)$.  Thus (cf. \cite{5,
Section~10}), the Ramsey property is equivalent to saying that
$\tau^{-1}(\uu,\uu)$ consists of exactly three points, namely the
ultrafilter on $\nn\times\nn$ mentioned above, a symmetrical one
``below diagonal,'' and an isomorphic copy of $\uu$ concentrated on
the diagonal.  (Hindman \cite{21} has shown that there are P-points
$\uu$ such that $\tau^{-1}(\uu,\uu)$ has the same cardinality,
$2^{2^{\aleph_0}}$, as $\bn$.)

P-points also have a Ramsey-like property \cite{2, 
Theorem 2.3}: If $\uu$ is a P-point and if $[\nn]^2$ is 
partitioned into two pieces, then there is $H\in\uu$ and 
there is a function $f:\nn\to\nn$ such that one 
piece of the partition contains all the two-element 
subsets $\{a,b\}$ of $H$ for which $f(a)<b$.

Although P-points and selective ultrafilters, like the 
ultrafilters discussed in the previous sections, have 
interesting combinatorial and topological properties, they 
are quite different in several respects, of which we list 
three.

First, since $\uu+\vv=\li\uu_n(n+\vv)$, while P-points 
are, in view of their topological description, never limit 
points of a countable set of other ultrafilters, it follows 
that no P-point can be of the form $\uu+\vv$.  In 
particular, no P-point can be recurrent.  So the family of 
P-points and, a fortiori, the subfamily of selective 
ultrafilters are disjoint from the families of ultrafilters 
studied in the preceding sections (recurrent, idempotent, 
etc.).

Second, the ultrafilters in the preceding sections are
proved to exist on the basis of the usual axioms of set 
theory (Zermelo-Fraenkel axioms, including the axiom of 
choice).  In contrast, the existence of P-points or of 
selective ultrafilters is independent of these axioms.  
More precisely, the continuum hypothesis (as well as 
weaker assumptions, like Martin's axiom) implies the 
existence of many selective ultrafilters and also many 
P-points that are not selective, but there are models of 
set theory with no selective ultrafilters \cite{25} and 
even with no P-points \cite{27, 28}.

Finally, where Ramsey ultrafilters have $\tau^{- 1}(\uu,\uu)$ as small
as possible, namely of size 3, the following theorem shows that
idempotent ultrafilters have it as large as possible.

\proclaim{Theorem  9}
If $\uu$ is an idempotent non-trivial ultrafilter on $\nn$, 
then there are $2^{2^{\aleph_0}}$ ultrafilters $\vv$ on 
$\nn\times\nn$ with $\tau(\vv)=(\uu,\uu)$.
\endproclaim

\demo{Proof}
Observe first that, for each natural number $n$, the set of 
multiples of $n$ is in $\uu$.  Indeed, as there are only 
finitely many congruence classes modulo $n$, any 
ultrafilter must contain one of them, so we can fix $j$ 
such that $0\leq j<n$ and $(\uu x)\,x\equiv j\pmod n$.  
Then 
$(\uu x)(\uu y)\,x+y\equiv 2j\pmod n$, so for 
idempotent $\uu$ it follows that $(\uu x)\,x\equiv 
2j\pmod n$.  But then $j\equiv 2j\pmod n$ and so $j=0$ 
as claimed.

In particular, for any $n$, $\uu$-almost all numbers $x$ 
are divisible by $2^n$ and therefore have 0's as the last 
$n$ digits in their binary expansions.

Using this, we can proceed as in Galvin's and Glazer's 
proof of Theorem 7 to find, in any set $C\in\uu$, a 
sequence $h_0,h_1,\dots$ such that 
\roster
\item All finite sums of distinct $h_i$'s are in $C$.
\item For each $i$, $h_{i+1}$ is divisible by a power of 2 
that is larger than $h_i$.
\endroster
Note the following consequence of \therosteritem2.  If 
$a$ and $b$ are each a sum of distinct $h_i$'s and if no 
$h_i$ occurs in both sums, then the 1's in the binary 
expansions of $a$ and $b$ occur in disjoint sets of 
positions.  We define the meshing number $m(a,b)$ to 
measure the amount of intermeshing between these 
disjoint sets; that is, $m(a,b)$ is the length $l$ of the 
longest sequence $s_1,\dots,s_l$ such that for all odd 
(resp. even) $i$, there is a 1 in position $s_i$ of the 
binary expansion of $a$ (resp. $b$).  It is clear that every 
integer $l\geq 2$ occurs as $m(a,b)$ with $a$ and $b$ 
sums of different $h_i$'s; just take the first $l$ of the 
$h_i$'s and let $a$ (resp. $b$) be the sum of the odd 
(resp. even) numbered ones.

In view of \therosteritem1, this means that each of the 
infinitely many sets $m^{-1}\{l\}$ meets each set of the 
form $C\times C$ for $C\in\uu$.  So each $m^{-1}\{l\}$ 
supports an ultrafilter containing all these $C\times C$ 
and therefore mapping to $(\uu,\uu)$ by $\tau$.  This 
proves that $\tau^{-1}(\uu,\uu)$ is infinite.  But every 
infinite closed subset of $\bn$ (or the homeomorphic 
$\beta(\nn\times\nn)$) has cardinality 
$2^{2^{\aleph_0}}$; see \cite{13, 18}.
\qed\enddemo

In spite of all these differences between P-points and 
selective ultrafilters on the one hand and recurrent and 
idempotent ultrafilters on the other, it is possible, using 
the continuum hypothesis (or Martin's axiom) to 
construct idempotent ultrafilters with strong connections 
to selective ultrafilters.  For example, one can arrange 
that an idempotent ultrafilter be mapped to a selective 
one by the map $\nn\to\nn$ that sends each natural 
number $a$ to the position of the rightmost (or the 
leftmost) 1 in its binary expansion.  For more information 
about such matters and for combinatorial applications, 
see \cite{7, 9}.

\enddocument